\documentclass[10pt]{article}
\usepackage{amsthm,amsmath,amssymb,amsfonts}
\usepackage{graphicx}

\usepackage[colorlinks=true,citecolor=black,linkcolor=black,urlcolor=blue]{hyperref}

\newcommand{\doi}[1]{\href{http://dx.doi.org/#1}{\texttt{doi:#1}}}
\newcommand{\arxiv}[1]{\href{http://arxiv.org/abs/#1}{\texttt{arXiv:#1}}}

\theoremstyle{plain}
\newtheorem{theorem}{Theorem}[section]
\newtheorem{lemma}[theorem]{Lemma}
\newtheorem{corollary}[theorem]{Corollary}
\newtheorem{proposition}[theorem]{Proposition}

\theoremstyle{definition}
\newtheorem{definition}[theorem]{Definition}
\newtheorem{example}[theorem]{Example}

\theoremstyle{remark}
\newtheorem{remark}[theorem]{Remark}

\def\p{{\mathsf{P}}}
\def\q{{\mathsf{Q}}}

\title{\bf Weighted $\mathsf{P}-$partitions enumerator}

\author{Marko Pe\v{s}ovi\'c\\
\small Faculty of Civile Engineering\\[-0.8ex]
\small University of Belgrade\\[-0.8ex]
\small\tt mpesovic@grf.bg.ac.rs\\
\and
Tanja Stojadinovi\'c\\
\small Faculty of Mathematics\\[-0.8ex]
\small University of Belgrade\\[-0.8ex]
\small\tt tanjas@matf.bg.ac.rs \and
Vladimir Gruji\'c\\
\small Faculty of Mathematics\\[-0.8ex]
\small University of Belgrade\\[-0.8ex]
\small\tt vgrujic@matf.bg.ac.rs}

\date{%\dateline{Jan 1, 2012}{Jan 2, 2012}\\
\small Mathematics Subject Classifications: 16T05, 52B05}

\begin{document}

\maketitle

\begin{abstract}
To an extended permutohedron we associate the weighted integer
points enumerator, whose principal specialization is the
$f$-polynomial. In the case of poset cones it refines Gessel's
$\mathsf{P}$-partitions enumerator. We show that this enumerator
is a quasisymmetric function obtained by universal morphism from
the Hopf algebra of posets.

\bigskip\noindent \textbf{Keywords}: generalized permutohedron, quasisymmetric function, poset, combinatorial Hopf algebra
\end{abstract}

\section{Introduction}

In the seminal paper of Aguiar, Bergeron and Sottile \cite{ABS}
the notion of combinatorial Hopf algebra was introduced. They
explained the ubiquity of quasisymmetric functions as generating
functions in enumerative combinatorics. More recently a geometric
meaning of quasisymmetric enumerators is attributed. It is based
on a class of convex polytopes called generalized permutohedra.
The integer points enumerator associated to a generalized
permutohedron is a quasisymmetric function. It was defined, and
studied in the case of matroid base polytopes, by Billera, Jia and
Reiner in \cite{BJR}, and in the case of nestohedra by Gruji\'{c}
in \cite{G}. More subtle generalization which takes into account
the face structure of a generalized permutohedron is introduced
and studied in \cite{GPS}. In this paper we consider the extended
generalized permutohedra and the special case of poset cones. We
prove that the integer points enumerator associated to a poset
cone coincides with the universal morphism from the Hopf algebra
of posets to quasisymmetric functions. The specialization for
$q=0$ is the Gessel enumerator of $\p$-partitions which attributes
the geometric meaning to this classical function.

In sections 2 we review necessary facts about quasisymmetric
functions and $\p$-partitions enumerators. In section 3 we
introduce the integer points enumerator $F_q(P)$ for an extended
generalized permutohedron $P$, which is a weighted quasisymmetric
function, in the same manner as in the case of generalized
permutohedra provided in \cite{GPS}. The parameter $q$ reflects
the rank function of the face lattice $L(P)$. To a poset $\p$ is
associated the poset cone $C(\mathsf{P})$, which is an extended
generalized permutohedron and our construction produces a weighted
quasisymmetric function $F_q(C(\mathsf{P}))$. In section 4 we
prove Theorem \ref{jakobitno}, the first main result of the paper,
which states that the weighted quasisymmetric function
$F_q(C(\mathsf{P}))$, constructed geometrically, has an algebraic
meaning as the universal morphism from a certain combinatorial
Hopf algebra of posets $\mathcal{P}$ to the Hopf algebra $QSym$ of
quasisymmetric functions. This result is analogous to the previous
results for simple graphs \cite{G1}, matroids and building sets
\cite{GPS}, and spreads their validity to the case of extended
generalized permutohedra. The main theorem is followed by various
examples, and statements about behavior of the enumerator of the
opposite poset and under the action of the antipode. In Theorem
\ref{druga} in section 5, it is shown that for a well labelled
poset $\p$ and $q=0$ our enumerator specializes to the classical
Gessel's $\p$-partitions enumerator. We also provide an example of
posets with the same $\p$-partitions enumerators but which are
distinguished by corresponding weighted quasisymmetric
enumerators.

\section{Quasisymmetric functions}

A \emph{composition} $\alpha$ of a positive integer $n$,
$\alpha\models n$, is an ordered list $(\alpha_1,\ldots,\alpha_k)$
of positive integers such that $\alpha_1+\cdots+\alpha_k=n$. The
\emph{monomial quasisymmetric function} $M_\alpha$ indexed by the
composition $\alpha$ is an element of the commutative algebra of
formal power series in the countable ordered set of variables
$\bold{x}=(x_1<x_2<x_3<\cdots)$ defined by

$$M_\alpha=\sum_{i_1<\cdots<i_k}x_{i_1}^{\alpha_1}\cdots x_{i_k}^{\alpha_k}.$$The \emph{algebra of quasisymmetric functions} $\mathcal{Q}Sym$,
spanned by $M_\alpha$ when $\alpha$ runs over all compositions, is
a subalgebra of the algebra of formal power series. The algebra
$\mathcal{Q}Sym$ is a graded, connected Hopf algebra (see
\cite{GR}, Proposition 5.8). The homogeneous component $QSym_n$ is
spanned by $\{M_\alpha\}_{\alpha\models n}$. Let
$\zeta_\mathcal{Q}:\mathcal{Q}Sym\rightarrow\bold{k}$ be a linear
multiplicative functional defined on the monomial basis by
$$\zeta_\mathcal{Q}(M_\alpha)=\begin{cases}1,&\text{if }\alpha=(n)\text{ for } n\in\mathbb{N},\\
0,&\text{otherwise.}\end{cases}$$ The Hopf algebra
$\mathcal{Q}Sym$ equipped with the character $\zeta_Q$ is the
terminal object in the category of combinatorial Hopf algebras.

\begin{theorem}[\cite{ABS}, Theorem 4.1]\label{GLAVNA}
For a combinatorial Hopf algebra $(\mathcal{H},\zeta)$ there is a
unique morphism of graded Hopf algebras
$\Psi:H\rightarrow\mathcal{Q}Sym$ such that
$$\Psi\circ\zeta_\mathcal{Q}=\zeta.$$

\noindent For a homogeneous element $h$ of degree $n$ the
coefficients $\zeta_\alpha(h)$,
$\alpha=(\alpha_1,\ldots,\alpha_k)\models n$ of $\Psi(h)$ in the
monomial basis  are given by
$$\zeta_\alpha(h)=\zeta^{\otimes k}\circ(p_{\alpha_1}\otimes\cdots\otimes p_{\alpha_k})\circ\Delta^{k-1}(h),$$

\noindent where $p_{i}$ is the projection of $\mathcal{H}$ on the
$i$-th homogeneous component $\mathcal{H}_i$ and $\Delta^{k-1}$ is
the $(k-1)-$fold coproduct map of $\mathcal{H}$.
\end{theorem}

For $F\in\mathcal{Q}Sym$ and $m\in\mathbb{N}$, let $\bold{ps}^1$
denotes the  \emph{principal specialization}
$$\bold{ps}^1(F)(m)=F(\underbrace{1,\ldots,1}_{m\text{}},0,0\ldots).$$
We have
$$\bold{ps}^1(M_\alpha)(m)=\binom{m}{k(\alpha)},$$
where $k(\alpha)$ is the number of parts of
$\alpha=(\alpha_1,\ldots,\alpha_k).$ Specially, for $m=-1$ we have
$$\bold{ps}^1(M_\alpha)(-1)=\binom{-1}{k(\alpha)}=(-1)^{k(\alpha)}.$$

For a composition $\alpha=(\alpha_1,\ldots,\alpha_k)\models n$ let
$D(\alpha)\subseteq[n-1]$ be a subset defined by
$D(\alpha)=\{\alpha_1,\alpha_1+\alpha_2,\ldots,
\alpha_1+\alpha_2+\cdots+\alpha_{k-1}\}$. We say that
$\alpha\models n$ \emph{refines} $\beta\models n$, and write
$\beta\preceq\alpha$, if $D(\beta)\subseteq D(\alpha).$\\Another
important  basis of $\mathcal{Q}Sym$ is the basis of
\emph{fundamental quasisymmetric functions} defined by
$$L_\alpha\;\;=\;\;\sum_{\alpha\preceq\beta}\;M_\beta.$$

\subsection{Quasisymmetric enumerator of $\mathsf{P}-$partitions}

A \emph{labelled poset} $\mathsf{P}$ is a poset on some finite
subset of positive integers. A \emph{$\mathsf{P}-$partition} is a
function $f:\mathsf{P}\rightarrow\mathbb{N}$ such that
\begin{enumerate}
\item[$\bullet$] $i<_\mathsf{P}j$ and $i<_\mathbb{Z}j$ implies
$f(i)\leq f(j)$, \item[$\bullet$] $i<_\mathsf{P}j$ and
$i>_\mathbb{Z}j$ implies $f(i)<f(j).$
\end{enumerate}

\begin{defn}\label{definicija}
A poset $\mathsf{P}$ is a \emph{well labelled poset} if
$i<_\mathsf{P}j$ implies $i>_\mathbb{Z}j.$ In that case
$\mathsf{P}-$partition is a function
$f:\mathsf{P}\rightarrow\mathbb{N}$ such that
$$i<_\mathsf{P}j\;\;\;\;\text{implies}\;\;\;\;f(i)<f(j).$$
Denote by $\mathcal{A}(\mathsf{P})$ the set of all
$\mathsf{P}-$partitions. Define \emph{the enumerator of
$\mathsf{P}-$partitions} by
$$F_\mathsf{P}(\bold{x})\;=\sum_{f\,\in\,\mathcal{A}(\mathsf{P})}x_{f(1)}x_{f(2)}\cdots x_{f(n)}.$$
\end{defn}

\begin{proposition}[\cite{GR}, Proposition 5.18]\label{t1}
For a totally ordered labelled poset
$\mathsf{P}=\{i_1<_\mathsf{P}i_2<_\mathsf{P}\cdots<_\mathsf{P}i_n\}$
the enumerator of $\mathsf{P}$-partitions is equal to the
fundamental quasisymmetric function
$$F_\mathsf{P}(\bold{x})=L_{\alpha(\mathsf{P})},$$
where $\alpha(\mathsf{P})\models n$ is a composition such that
$D(\alpha(\mathsf{P}))=\{j:i_j>_\mathbb{Z}i_{j+1}\}.$
\end{proposition}

\begin{theorem}[\cite{GR}, Theorem 5.19]\label{t2}
For a labelled poset $\p$,
$$F_\p(\bold{x})=\sum_{\mathsf{l}\,\in\,\mathcal{L}(\p)}F_\mathsf{l}(\bold{x}),$$
where the sum is over the set $\mathcal{L}(\p)$ of all linear
extensions $\mathsf{l}$ of $\p$.
\end{theorem}

\begin{example}
Figure 1 presents two posets with their enumerators of
$\p$-partitions. Note that the poset $\p_1$ is not well labelled,
while the poset $\p_2$ is.
\begin{figure}[h!]
\begin{center}
\includegraphics[width=50mm]{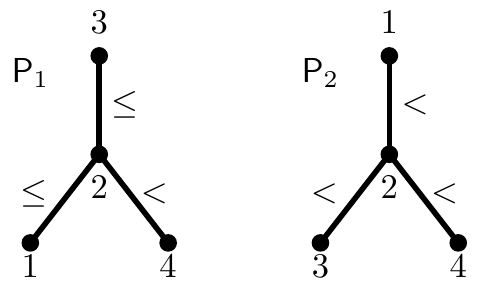}
\caption{$F_{\p_1}(\bold{x})=L_{(1,3)}+L_{(2,2)}$ and
$F_{\p_2}(\bold{x})=L_{(2,1,1)}+L_{(1,1,1,1)}$}
\end{center}
\end{figure}
\end{example}

\section{Extended generalized permutohedra}

A \emph{standard $(n-1)-$dimensional permutohedron} $Pe^{n-1}$ is
the convex hull of the orbit of a point
$(a_1,a_2,\ldots,a_n)\in\mathbb{R}^n$ with increasing coordinates
$a_1<a_2<\cdots<a_n$
$$Pe^{n-1}=\mathrm{conv}\{(a_{\sigma(1)},a_{\sigma(2)},\ldots,a_{\sigma(n)})\,\,:\,\,\sigma\in\mathfrak{S}_n\},$$
\noindent where $\mathfrak{S}_n$ is the permutation group of
$[n]$.

\begin{proposition}[\cite{P}, Proposition 2.6]
The $d-$dimensional faces of $Pe^{n-1}$ are in one-to-one
correspondence with set compositions
$\mathcal{C}=C_1|C_2|\cdots|C_{n-d}$ of $[n].$ The face
corresponding to the set composition $\mathcal{C}$ is given by the
$n-d$ linear equations
$$\sum_{i\in C_1\cup C_2\cup\cdots\cup C_k}x_i\;\;=\;\;a_1+a_2+\cdots+a_{|C_1\cup C_2\cup\cdots\cup C_k|},$$
for $1\leq k\leq n-d$.
\end{proposition}

A set composition $\mathcal{C}=C_1|C_2|\cdots|C_{n-d}$ defines the
\emph{flag} $\mathcal{F}$ of subsets
$$\mathcal{F}\;:\;\emptyset=F_0\subset F_1\subset\cdots\subset F_{n-d-1}\subset F_{n-d}=[n],$$
of the length $|\mathcal{F}|=n-d$, where $F_i=C_1\cup
C_2\cup\cdots\cup C_i$, for $1\leq i\leq n-d$. There is an obvious
order reversing one$-$to$-$one correspondence between the face
lattice of the permutohedron $L(Pe^{n-1})$ and the lattice of
flags of subsets of the set $[n].$ Using this correspondence we
will label faces of the standard permutohedron $Pe^{n-1}$ by flags
of subsets of $[n]$. We have
$$\dim(\mathcal{F})=n-|\mathcal{F}|.$$

The normal fan $\mathcal{N}(Pe^{n-1})$ of the standard
permutohedron $Pe^{n-1}$ is the fan of the \emph{braid
arrangement}, given by hyperplanes $\{x_i=x_j\}_{1\leq i<j\leq
n}$, in $\mathbb{R}^n$. The cones of the braid arrangement fan are
called \emph{braid cones.} The braid cone $C_\mathcal{F}$ at the
face $\mathcal{F}$ is determined by
\begin{enumerate}
\item[$\bullet$] $x_p=x_q$ if $p,q\in F_{i+1}\setminus F_i$, for
some $0\leq i\leq |\mathcal{F}|-1,$ \item[$\bullet$] $x_p\leq x_q$
if $p\in F_i\setminus F_{i-1}$ and $q\in F_{i+1}\setminus F_i$,
for some $1\leq i\leq |\mathcal{F}|-1$.
\end{enumerate}
We have $\dim(C_\mathcal{F})=|\mathcal{F}|$. The fan
$\mathcal{N}_1$ is a \emph{refinement} of $\mathcal{N}_2$ (or
$\mathcal{N}_2$ is a \emph{coarsement} of $\mathcal{N}_1$) if
every cone of $\mathcal{N}_1$ is  contained in a cone in
$\mathcal{N}_2$ (or if every cone in $\mathcal{N}_2$ is a union of
cones of $\mathcal{N}_1$).

\begin{defn}
A convex polytope $Q$ is an $(n-1)$-dimensional \emph{generalized
permutohedron} if the braid arrangement fan
$\mathcal{N}(Pe^{n-1})$ refines the normal fan $\mathcal{N}(Q)$.
\end{defn}

For a generalized permutohedron $Q$ there is a map
$\pi_Q:L(Pe^{n-1})\rightarrow L(Q)$ between face latices,
determined by $\pi_Q(\mathcal{F})=G$ if and only if the relative
interior of the braid cone $C^\circ_\mathcal{F}$ is contained in
the relative interior $C^\circ_G$ of the normal cone $C_G$ at the
face $G\in L(Q)$. We say that the flag $\mathcal{F}$ is {\it
normal} to the face $G$. Denote by
$\mathsf{F}(G)=\{\mathcal{F}:C^\circ_\mathcal{F}\subseteq
C^\circ_{G}\}$ the set of normal flags to a face $G$. By
Proposition 2.3 in \cite{GPS}, we have
\begin{equation}\label{sssss}
\sum_{\mathcal{F}\in
\mathsf{F}(G)}(-1)^{|\mathcal{F}|}\;\,=\,\;(-1)^{n-1-\dim(G)}.
\end{equation}

\begin{defn}
An \emph{extended generalized permutohedron} $P$ is a polyhedron
whose normal fan $\mathcal{N}(P)$ is a coarsening of a subfan of
the braid arrangement fan.
\end{defn}

\subsection{Quasisymmetric enumerator $F_q(P)$}

A vector
$\omega=(\omega_1,\omega_2,\ldots,\omega_n)\in\mathbb{Z}^n_+$
defines the \emph{weight function}
$\omega^\ast:\mathbb{R}^n\rightarrow\mathbb{R}$  by
$\omega^\ast(x):=\left<\omega,x\right>$ where $\left<\,,\right>$
is the standard scalar product in $\mathbb{R}^n$. The weight
function $\omega^\ast$ on $Pe^{n-1}$ is maximized along a unique
face $\mathcal{F}_\omega$ of $Pe^{n-1}$ determined by the
condition that the vector $\omega$ lies in the relative interior
of its normal cone $\omega\in C^\circ_{\mathcal{F}_\omega}$. The
flag $\mathcal{F}_\omega$ satisfies the following conditions
\begin{enumerate}
\item[$\bullet$] $\omega$ is constant on $F_i\setminus F_{i-1}$,
for $1\leq i\leq k$, \item[$\bullet$] $\omega|_{F_i\setminus
F_{i-1}}<\,\omega|_{F_{i+1}\setminus F_i},$ for $1\leq i\leq k-1$,
\end{enumerate}
where $k$ is the length of $\mathcal{F}_\omega$. Let
$M_\mathcal{F}$ be the \emph{enumerator} of positive integer
vectors $\omega\in\mathbb{Z}^n_+$ in relative interior of the
corresponding braid cone $C_\mathcal{F}$
\begin{equation}\label{jedinicica}
M_\mathcal{F}\;\;=\sum_{\omega\,\in\,\mathbb{Z}^n_+\,\cap\,\,C^\circ_\mathcal{F}}x_{\omega_1}x_{\omega_2}\cdots
x_{\omega_n}.
\end{equation}
Note that $\omega\in C^\circ_\mathcal{F}$ if and only if
$\mathcal{F}_\omega=\mathcal{F}$. The enumerator $M_\mathcal{F}$
is a monomial quasisymmetric function indexed by composition
$$\mathrm{type}(\mathcal{F})=(|F_1\setminus F_0|, |F_2\setminus F_1|,\ldots, |F_k\setminus F_{k-1}|).$$
A \emph{weight function} $\omega^\ast$ on an extended generalized
permutohedron $P$ is maximized along a unique face $G_\omega$ of
$P$ which is determined by the condition $\omega\in
C^\circ_{G_\omega}$. For a face $G$ of $P$ we define a
\emph{quasisymmetric enumerator} $F_q(G)$ by
$$F_q(G)=q^{\dim(G)}\sum_{\omega\in\mathbb{Z}^n_+\,\cap\,C^\circ_G}x_{\omega_1}x_{\omega_2}\cdots x_{\omega_n}.$$
By (\ref{jedinicica}) it follows immediately that
$$F_q(G)=q^{\dim(G)}\sum_{\mathcal{F}\in\,\mathsf{F}(G)}M_\mathcal{F}.$$

\begin{defn}\label{enumerator}
For an extended generalized permutohedron $P$ with the face
lattice $L(P)$ let
$$F_q(P)=\sum_{G\in L(P)}F_q(G).$$
\end{defn}

Let
$f(P,q):=f_0\,+\,f_1q\,+\,f_2q^2\,+\,\cdots\,+\,f_{n-1}q^{n-1}$ be
the \emph{$f-$polynomial} of $(n-1)-$dimensional polyhedron $P$.
The coefficient $f_i$ is the number of $i-$dimensional faces of
polyhedron $P$. The following proposition describes the
$f-$polynomial of extended generalized permutohedron $P$ in terms
of the weighted quasisymmetric enumerator $F_q(P).$

\begin{proposition}\label{f-f-f}
The $f-$polynomial of an $(n-1)-$dimensional extended generalized
permutohedron $P$ is determined by the principal specialization
\begin{equation*}
f(P,q)\,=\,(-1)^{n-1}\bold{ps}^1(F_{-q}(P))(-1).
\end{equation*}
\end{proposition}
\noindent\textit{Proof.} It follows from (\ref{sssss}) that

$$\bold{ps}^1(F_{-q}(G))(-1)=(-1)^{n-1}q^{\dim(G)}.$$
Therefore

\[\bold{ps}^1(F_{-q}(P))(-1)=(-1)^{n-1}\displaystyle\sum_{G\in
L(P)}q^{\dim(G)}.\]

\subsection{Poset cone}\label{sekcija}

For a poset $\p$ on $[n]$ denote by $\p|_S$ the \emph{restriction}
of $\p$ to $S\subseteq[n]$. A subset $S\subseteq[n]$ is called
\emph{ideal} of $\p$, denoted by $S\lhd \p$, if no element of
$[n]\setminus S$ is less than an element of $S$.

\begin{defn}
The \emph{poset cone} of a poset $\p$ on $[n]$ is an extended
generalized permutohedron given by
$$C(\mathsf{P})=\mathrm{cone}\{e_i-e_j:i<_\p j\},$$ where
$e_1,\,e_2,\,\ldots,\,e_n$ are the standard basis vectors in
$\mathbb{R}^n$.
\end{defn}

\begin{proposition}[\cite{HM} Proposition 15.1]
The generating rays of $C(\mathsf{P})$ are determined by the
vectors $e_i-e_j$ corresponding to the cover relations
$i\lessdot_\p j$ in $\p$.
\end{proposition}

It is shown in \cite{HM} that the poset cone is described by
$$C(\mathsf{P})=\left\{(x_1,\ldots,x_n):\sum_{i=1}^nx_i=0\;\;\text{and}\;\;\sum_{s\in S}x_s\geq0\;\;\text{for all}\;\;S\lhd \p\right\}.$$
Therefore
\begin{equation}\label{dimenzija}
\dim(C(\mathsf{P}))=n-c(\p),
\end{equation}
where $c(\p)$ is the number of connected components of the poset $\p$.\\

Let $\mathsf{C}\,:\,i_1,\,i_2,\,\ldots,\,i_n$  be a cyclic
sequence  of elements of $\p$ where every consecutive pair is
comparable in $\p.$ We say that $\mathsf{C}$ is a \emph{circuit}
of $\p$. Circuits consist of \emph{up$-$edges}  $i_j<_\p i_{j+1}$
and \emph{down$-$edges} $i_j >_\p i_{j+1}$.

\begin{defn}
A subposet $\mathsf{Q}$ of a poset $\p$ on $[n]$ is
\emph{positive} if for every circuit $\mathsf{C}$ all down$-$edges
of $\mathsf{C}$ are in $\mathsf{Q}$ if and only if all up$-$edges
of $\mathsf{C}$ are in $\mathsf{Q}$. Let $\mathrm{Pos}(\p)$ be the
set of all positive subposets of $\p$.
\end{defn}

The faces of the poset cone $C(\mathsf{P})$ are characterized by
the following lemma.
\begin{lemma}[\cite{HM}, Lemma 15.3]\label{strane}
Let $\p$ be a poset on $[n]$. The faces of the poset cone
$C(\mathsf{P})\subseteq \mathbb{R}^n$ are precisely the poset
cones $C(\mathsf{Q})$ as $\mathsf{Q}$ ranges over positive
subposets of $\p$.
\end{lemma}

By Definition \ref{enumerator} and Lemma \ref{strane}, it follows
that the weighted quasisymmetric enumerator $F_q(C(\mathsf{P}))$
for a poset $\p$ can be expressed as
\begin{equation}\label{enumerator_p}
F_q(C(\mathsf{P}))=\sum_{\mathsf{Q}\in\mathrm{Pos}(\p)}
\sum_{\mathcal{F}\in\mathsf{F}(C(\mathsf{Q}))}q^{\dim(C(\mathsf{Q}))}M_\mathcal{F}.
\end{equation}

\begin{proposition}\label{stav_1}
A vector $\omega=(\omega_1,\ldots,\omega_n)\in\mathbb{Z}^n_+$ lies
in $\mathcal{N}(C(\mathsf{P}))$ if and only if
$\,\omega_i\leq\omega_j$ for any $i\leq_\p j$.
\end{proposition}
\noindent\textit{Proof.} A vector $\omega\in\mathbb{Z}^n_+$ lies
in $\mathcal{N}(C(\mathsf{P}))$ if and only if $\omega^\ast$ is
maximized over $C(\mathsf{P})$ at some its face. The weight
function $\omega^\ast$ on the ray generated by $e_i-e_j$ is given
by
$$\omega^\ast(t(e_i-e_j))=\left<\omega,t(e_i-e_j)\right>=t(\omega_i-\omega_j), t\geq0.$$
Consequently, we have:
\begin{enumerate}
\item If $\omega_i>\omega_j$, then $t(\omega_i-\omega_j)\geq0$ for
$t\geq0$ and maximum of $\omega^{\ast}$ is not achieved along the
ray generated by $e_i-e_j$.

\item If $\omega_i<\omega_j$, then $t(\omega_i-\omega_j)\leq0$ for
$t\geq0$ and maximum along the ray generated by $e_i-e_j$ is
achieved at the vertex of $C(\mathsf{P})$.

\item If $\omega_i=\omega_j$, then $t(\omega_i-\omega_j)=0$ for
$t\geq0$ and $\omega^\ast$ is maximized along the ray generated by
$e_i-e_j$.\qed
\end{enumerate}

\begin{corollary}\label{posledica_za_dokaz}
A flag $\mathcal{F}$ is normal to the face $C(\q)$ for some
positive subposet $\q\in\mathrm{Pos}(\p)$, i.e.
$\mathcal{F}\in\mathsf{F}(C(\q))$ if and only if for each
$\omega\in C^\circ_\mathcal{F}$ holds that
\begin{enumerate}
\item $\omega_i=\omega_j$ for all $i\leq_\q j$, \item
$\omega_i<\omega_j$ for all $i,j$ which are incomparable in $\q$
and $i\leq_\p j$.
\end{enumerate}
\end{corollary}

The following examples illustrates the concepts introduced in
connection with poset cones.

\begin{example}

Let $\p$ be a poset on $[4]$ defined by covering relations
$1,2<3,4$. The generated rays of the poset cone $C(\p)$ are
determined by vectors $e_1-e_3,e_1-e_4,e_2-e_3,e_2-e_4$. The
normal fan $\mathcal{N}(C(\p))$ is described by inequalities
$\omega_1,\omega_2\leq\omega_3,\omega_4$, see Figure 2. The list
of circuits of $\p$ is the following

\[1<4>2<3, 1<3>2<4, 4>2<3>1, 3>2<4>1,\] \[2<3>1<4, 2<4>1<3, 3>1<4>2,
4>1<3>2.\] Its facets are determined by positive subposets given
by \[1<3,4; 2<3,4; 1,2<3; 1,2<4.\] Direct calculation gives
$F_q(C(\p))=q^{3}M_{(4)}+2q^{2}(M_{(1,3)}+M_{(3,1)})+4qM_{(1,2,1)}+M_{(2,2)}+2M_{(1,1,2)}+2M_{(2,1,1)}+4M_{(1,1,1,1)}.$

\begin{figure}[h!]
\begin{center}
\includegraphics[width=60mm]{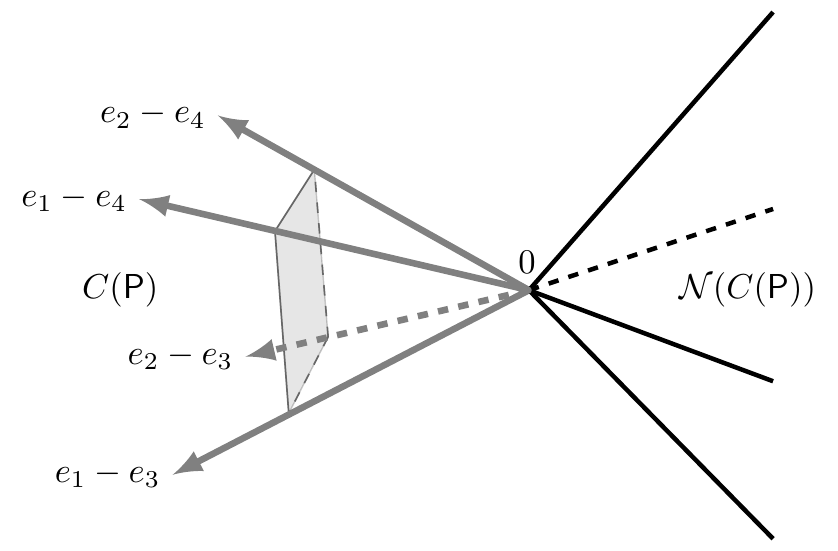}
\caption{Poset cone and its normal fan}
\end{center}
\end{figure}

\end{example}

\section{Hopf algebra $\mathcal{P}$}

The set of isomorphism classes of finite posets linearly generates
the $\bold{k}-$vector space $\mathcal{P}$
$$\mathcal{P}=\bigoplus_{n\geq0}\mathcal{P}_n,$$
where $\mathcal{P}_n$ is the homogeneous component of degree $n$
spanned by posets on $[n].$ The space $\mathcal{P}$ is a graded,
commutative and non-cocommutative Hopf algebra (see \cite{ABS},
Example 2.3) with the \emph{multiplication}
$$[\p_1]\cdot [\p_2] = [\p_1\sqcup \p_2],$$
where $\p_1\sqcup \p_2$ is disjoint union of posets, and the
\emph{comultiplication}
$$\Delta([\p])=\sum_{S\,\lhd\,\p}[\p|_S]\otimes [\p|_{[n]\setminus S}].$$

\begin{defn}
A flag $\mathcal{F}:\emptyset=F_0\subset F_1\subset\cdots\subset
F_k=[n]$ is a \emph{flag of ideals} of а poset $\p$, denoted by
$\mathcal{F}\lhd \p$, if $F_i\lhd \p,$ for all $0< i\leq k.$ Let
$$\mathfrak{F}(\p)=\{\mathcal{F}:\mathcal{F}\lhd \p\}.$$
\end{defn}
There is a map $\mathfrak{F}(\p)\rightarrow\mathcal{P}_n$ given by
$$\mathcal{F}\mapsto\p/\mathcal{F}=\prod_{i=1}^k\p|_{F_i\setminus F_{i-1}}.$$
A poset $\p/\mathcal{F}$ is a weak subposet of $\p$, i.e. if
$i\leq j$ in $\p/\mathcal{F}$ then $i\leq j$ in $\p$. Let
$\p/\mathcal{F}=\p_1\sqcup \p_2\sqcup\cdots\sqcup \p_m$ be the
decomposition into connected components. We say that $\p_u$ and
$\p_v$ are \emph{incomparable} if elements from $\p_u$ and $\p_v$
are mutually incomparable in a poset $\p$. Otherwise, $\p_u$ is
\emph{smaller} than $\p_v$ if no element of $\p_v$ is less than an
element of $\p_u$.

\subsection{Quasisymmetric enumerator $F_q(C(\mathsf{P}))$}

We extend the basic field $\bold{k}$ into the field of rational
functions $\bold{k}(q)$ and define the character
$\zeta_q:\mathcal{P}\rightarrow\bold{k}(q)$ with
$$\zeta_q([\p])=q^{\mathrm{rk}(\p)}\;\;\;\;\text{where}\;\;\;\;\mathrm{rk}(\p)=n-c(\p).$$

Let
$\Psi_q:(\mathcal{P},\zeta_q)\rightarrow(\mathcal{Q}Sym,\zeta_\mathcal{Q})$
be the unique morphism of combinatorial Hopf algebras over
$\bold{k}(q)$, given by Theorem \ref{GLAVNA} with
\begin{equation}\label{tri}
\Psi_q([\p])=\sum_{\alpha\models n}(\zeta_q)_\alpha(\p)M_\alpha.
\end{equation}
The coefficient corresponding to a composition
$\alpha=(\alpha_1,\alpha_2,\ldots,\alpha_k)\models n$ is
determined by
$$(\zeta_q)_\alpha(\p)=\sum_{\substack{\mathcal{F}\in\mathfrak{F}(\p)\\\mathrm{type}(\mathcal{F})=\alpha}}
\prod_{j=1}^kq^{\mathrm{rk}(\p|_{F_j\setminus
F_{j-1}})}=\sum_{\substack{\mathcal{F}\in\mathfrak{F}(\p)\\\mathrm{type}(\mathcal{F})=\alpha}}q^{\mathrm{rk}_\p(\mathcal{F})},$$
where
$$\mathrm{rk}_\p(\mathcal{F})=\sum_{j=1}^k\mathrm{rk}(\p|_{F_j\setminus F_{j-1}})=n-\sum_{j=1}^kc(\p|_{F_j\setminus F_{j-1}}).$$
Thus, the equation (\ref{tri}) can be expressed as
\begin{equation*}
\Psi_q([\p])=\sum_{\mathcal{F}\in\mathfrak{F}(\p)}q^{\mathrm{rk}_\p(\mathcal{F})}M_\mathcal{F}.
\end{equation*}

To a poset $\p$ are associated two weighted quasisymmetric
functions, $F_q(C(\mathsf{P}))$ with a geometric meaning based on
the combinatorics of the poset cone $C(\mathsf{P})$ and
$\Psi_q([\p])$ with an algebraic meaning based on the Hopf algebra
structure on finite posets. The following theorem shows that these
two functions are equal.

\begin{theorem}\label{jakobitno}
For a poset $\p$ the quasisymmetric enumerator function
$F_q(C(\mathsf{P}))$ associated to a poset cone $C(\mathsf{P})$
coincides with the value of the universal morphism at $[\p]$ from
the combinatorial Hopf algebra of posets $\mathcal{P}$ to
$\mathcal{Q}Sym$
$$F_q(C(\mathsf{P}))=\Psi_q([\p]).$$
\end{theorem}
\noindent\textit{Proof.} The weighted enumerator
$F_q(C(\mathsf{P}))$ is described by (\ref{enumerator_p}). We need
to show that
\[\mathfrak{F}(\p)=\{\mathcal{F}:\mathcal{F}\in\mathsf{F}(C(\q))\text{ for some }\q\in\mathrm{Pos}(\p)\},\]
and $\mathrm{rk}_\p(\mathcal{F})=\dim(C(\q))$, for $\mathcal{F}\in\mathsf{F}(C(\q))$. Let $\p$ be a poset on the set $[n]$.\\

$\supseteq:$ Suppose $\mathcal{F}\in\mathsf{F}(C(\q))$ for some
$\q\in\mathrm{Pos}(\p)$. If $\q=\q_1\sqcup\ldots\sqcup \q_m$ is
the decomposition into connected components, by (\ref{dimenzija})
we have $\dim(C(\q))=n-m$. For $\omega\in C^\circ_\mathcal{F}$, by
Corollary \ref{posledica_za_dokaz}, we deduce the following facts:
\begin{enumerate}
\item[$\bullet$]
 $\omega$ is constant on each $\q_u$, hence $\q_u\subseteq F_i\setminus F_{i-1}$ for some
 $i=1,\ldots,|\mathcal{F}|$.

\item[$\bullet$] if $\q_u, \q_v\subset F_{i}\setminus F_{i-1}$ for
some $i=1,\ldots,|\mathcal{F}|$, then $\q_u$ and $\q_v$ are
incomparable. Otherwise, if $i\leq_\p j$ for some
$i\in\q_u,j\in\q_v$, then $\omega^\ast$ is maximized along the ray
generated by $e_i-e_j$, contrary to
$\mathcal{F}\in\mathsf{F}(C(\q))$.

\item[$\bullet$] if $\q_u\subseteq F_i\setminus F_{i-1}$,
$\q_v\subseteq F_j\setminus F_{j-1}$ and $\q_u$ is smaller than
$\q_v$, then $i<_\mathbb{Z}j$. Otherwise, if $i>_\mathbb{Z}j$,
then $\omega_i>\omega_j$ for some $i\in \q_u$, $j\in \q_v$ such
that $i<_\p j$, hence $\omega^\ast$ does not reach the maximum.
\end{enumerate}
We conclude that $\mathcal{F}\in\mathfrak{F}(\p)$ and
$\p/\mathcal{F}=\q$, consequently
$\mathrm{rk}_\p(\mathcal{F})=n-m=\dim(C(\q))$.

$\subseteq:$ Let $\mathcal{F}\lhd\p$ be a flag of ideals of the
poset $\mathsf{P}$ and $\mathsf{C}:i_1,i_2,\ldots,i_n$ be a
circuit with all down-edges belonging to $\mathsf{P}/\mathcal{F}$.
The level function of elements of $\p$ according to the flag
$\mathcal{F}$, defined by $l(i)=\min\{a\mid i\in F_a\},
i\in\mathsf{P}$, is nondecreasing along the circuit $\mathsf{C}$.
It means that $l$ is constant on $\mathsf{C}$, i.e.
$\mathsf{C}\subset F_a\setminus F_{a-1}$, for some $1\leq
a\leq|\mathcal{F}|$. Particulary, all up-edges of the circuit
$\mathsf{C}$ are in $\mathsf{P}/\mathcal{F}$. The same is true for
circuits with all up-edges in $\mathsf{P}/\mathcal{F}$ and we
conclude that $\mathsf{P}/\mathcal{F}$ is a positive subposet of
$\p$. Corollary \ref{posledica_za_dokaz} gives
$\mathcal{F}\in\mathsf{F}(C(\mathsf{P}/\mathcal{F}))$ and the
proof is finished by the obvious identity
$\mathrm{rk}_\mathsf{P}(\mathcal{F})=\dim
C(\mathsf{P}/\mathcal{F})$. \qed

By Theorem \ref{jakobitno} and Proposition \ref{f-f-f} we obtain
the following expression for the $f$-polynomial.

\begin{corollary}\label{fpsi}
Let $C(\mathsf{P})$ be the poset cone associated to a poset $\p$
on the ground set $[n].$ The $f-$polynomial of $C(\mathsf{P})$ is
given by
$$f(C(\mathsf{P}),q)=(-1)^{n-1}\bold{ps}^1(\Psi_{-q}([\p]))(-1).$$
\end{corollary}

\begin{example}
Let $\mathsf{st}_n$ be the poset on $[n]$ with covering relations
$i\lessdot n$, for $1\leq i\leq n-1$.

%\begin{figure}[h!]
%\begin{center}
%\includegraphics[width=60mm]{slika_2.pdf}
%\caption{The posets $\mathsf{st}_2, \mathsf{st}_3, \mathsf{st}_4$}
%\end{center}
%\end{figure}

Let $\mathcal{F}:\emptyset=F_0\subset F_1\subset\cdots\subset
F_k=[n]$ be a flag of ideals of the poset $\mathsf{st}_n$. Then
$n\in F_k$ and two different elements of $[n]$ are in the same
connected component in $\mathsf{st}_n/\mathcal{F}$ if and only if
the both are in $F_k\setminus F_{k-1}$. Therefore

$$\mathrm{rk}_{\mathsf{st}_n}(\mathcal{F})=|F_k\setminus F_{k-1}|-1.$$
The number of flags corresponding to the faces of
$C(\mathsf{st}_n)$ of dimension $i$ is equal to $\binom{n-1}{i}$,
which implies
\begin{equation*}
F_q(C(\mathsf{st}_n))=\sum_{i=0}^{n-1}\binom{n-1}{i}\left(M_{(1)}^{n-1-i}\right)_{i+1}q^i,
\end{equation*}
where for $F\in\mathcal{Q}Sym$, $F\mapsto (F)_i$ is the linear
extension of the map given on monomial basis by $M_\alpha\mapsto
M_{(\alpha,i)}$. By Proposition \ref{fpsi}, the corresponding
$f-$polynomial is equal to

\[f(C(\mathsf{st}_n),q)=\displaystyle\sum_{i=0}^{n-1}\binom{n-1}{i}q^i=(1+q)^{n-1}.\]
\end{example}

\begin{example}
Let $\mathsf{l}_n$ be a linear poset on $[n]$ and
$\mathcal{F}\in\mathfrak{F}(\mathsf{l}_n)$. For $1\leq
i\leq|\mathcal{F}|$ all components of $F_i\setminus F_{i-1}$ are
connected, so
$$c(F_{i}\setminus F_{i-1})=1\;\;\text{and}\;\;\mathrm{rk}(\mathsf{l}|_{F_{i+1}\setminus F_i})=|F_{i+1}\setminus F_i|-1.$$
It implies that
$\mathrm{rk}_{\,\mathsf{l}_n}(\mathcal{F})=n-|\mathcal{F}|$. We
have
$$F_q(C(\mathsf{l}_n))=\sum_{i=0}^{n-1}\left(\sum_{\alpha\,:\,k(\alpha)=n-i}M_\alpha\right)q^i.$$
Since $|\{\alpha\models n\,:\,k(\alpha)=n-i\}|=\binom{n-1}{i},$
the corresponding $f$-polynomial is equal to

\[f(C(\mathsf{l}_n),q)=\displaystyle\sum_{i=0}^{n-1}\binom{n-1}{i}q^i=(1+q)^{n-1}.\]
\end{example}
We obtain that $f(C(\mathsf{st}_n),q)=f(C(\mathsf{l}_n),q)$.
Actually, it is a consequence of a more general fact.

\begin{proposition}
Let $\mathsf{P}$ be a poset on $[n]$ whose Hasse diagram is a
tree. Then
$$f(C(\mathsf{P}),q)=(1+q)^{n-1}.$$
\end{proposition}
\noindent\textit{Proof.} The poset $\mathsf{P}$ has $n-1$ covering
relations. Proposition \ref{strane} implies that the generating
rays of $C(\mathsf{P})$ are $n-1$ linearly independent vectors
$e_i-e_j$ where $j\lessdot i$. Hence, the $k-$face of
$C(\mathsf{P})$ is generated by $k$ generating rays, so
$f_k(C(\mathsf{P}))=\binom{n-1}{k}$.\qed

\begin{example}\label{bipart}
Let $\mathsf{K}_{m,n}$ be the poset on the set $[m+n]$ such that
for all $i\in[m]$ and $j\in[m+n]\setminus[m]$ hold $i\lessdot j.$
The Hasse diagram of $\mathsf{K}_{m,n}$ is the complete bipartite
graph $K_{m,n}$.

%\begin{figure}[h!]
%\begin{center}
%\includegraphics[width=60mm]{slika_3.pdf}
%\caption{Posets $\mathsf{K}_{3,3}$ and $\mathsf{K}_{3,2}$}
%\end{center}
%\end{figure}

We have
\begin{equation*}\begin{split}
F_q(C(\mathsf{K}_{m,n}))& =\left(M_{(1)}^m\right)\circ
\left(M_{(1)}^n\right)+\\
&+\sum_{k=1}^{m+n-1}\;\;\;q^k\sum_{t_1+t_2=k+1}\binom{m}{t_1}\binom{n}{t_2}\left(M_{(1)}^{m-t_1}\circ
M_{(k+1)}\circ M_{(1)}^{n-t_2}\right),
\end{split}\end{equation*}
where $1\leq t_1\leq m$, $1\leq t_2\leq n$ and $\circ$ is the
\emph{concatenation product} defined on monomial basis by
$M_\alpha\circ M_\beta:=M_{\alpha\,\cdot\,\beta}$. Note that this
includes the case $\mathsf{st}_n=\mathsf{K}_{n-1,1}$. The
principal specialization evaluated at $-1$ gives
\[
f(C(\mathsf{K}_{m,n},q)=1+\sum_{k=1}^{m+n-1}\;\;\;q^k\sum_{t_1+t_2=k+1}\binom{m}{t_1}\binom{n}{t_2}.
\]
\end{example}

\subsubsection{Opposite poset}

For a flag $\mathcal{F}:\emptyset=F_0\subset
F_1\subset\cdots\subset F_k=[n]$, the \emph{opposite flag}
$\mathcal{F}^{op}$ is defined by
$$\mathcal{F}^{op}:\emptyset\subset[n]\setminus F_{k-1}\subset[n]\setminus F_{k-2}\subset\cdots\subset[n]\setminus F_1\subset[n].$$
The normal cone corresponding to the opposite flag is the opposite
cone
$$C_{\mathcal{F}^{op}}=-C_{\mathcal{F}}\,\,\text{and}\,\,\mathrm{type}(\mathcal{F}^{op})=\mathrm{rev}(\mathrm{type}(\mathcal{F})),
$$
where
$\mathrm{rev}(\alpha_1,\alpha_2,\ldots,\alpha_k)=(\alpha_k,\alpha_{k-1},\ldots,\alpha_1).$
The \emph{opposite poset} $\p^{op}$ to a poset $\p$ is the poset
on the same set  such that $i<_{\p^{op}}j$ if and only if
$j<_{\p}i.$

\begin{proposition}
Let $\p$ be a poset on $[n]$, the quasisymmetric enumerator
function corresponding to the opposite poset $\p^{op}$ is
determined by
$$F_q(C(\p^{op}))\;\,=\,\;\mathrm{rev}(F_q(C(\p)),$$
%=\sum_{\mathcal{F}\in\,\mathfrak{F}(\p)}q^{\mathrm{rk}_\p(\mathcal{F})}M_{\mathcal{F}^{op}}.$$
where for $F\in\mathcal{Q}Sym$, $F\mapsto\mathrm{rev}(F)$ is the
linear extension of the map given on monomial basis by
$M_\alpha\mapsto M_{\mathrm{rev}(\alpha)}$.
\end{proposition}
\noindent\textit{Proof.} Note that
$\mathcal{F}\in\mathfrak{F}(\p)$ if and only if
$\mathcal{F}^{op}\in\mathfrak{F}(\p^{op}).$ The statement follows
from
$\mathrm{rk}_\p(\mathcal{F})=\mathrm{rk}_{\p^{op}}(\mathcal{F}^{op}).$\qed

\subsubsection{Action of the antipode on $F_q(C(\mathsf{P}))$}

The geometric interpretation of the action of the antipode on
$\mathcal{Q}Sym$ is given by the following lemma.

\begin{lemma}[\cite{GPS}, Lemma 4.5]
The antipode $S$ on the monomial quasisymmetric function
$M_\mathcal{F}$ associated to a flag $\mathcal{F}$ acts by
$$S(M_\mathcal{F})\;=\;(-1)^{|\mathcal{F}|+1}\sum_{\mathcal{G}\,\preceq\,\mathcal{F}^{op}} M_\mathcal{G},$$
where $\mathcal{G}\,\preceq\,\mathcal{F}^{op}$ if and only if
$\mathcal{F}^{op}\subseteq\mathcal{G}$ as faces of $Pe^{n-1}$.
\end{lemma}

The next theorem extends the similar statement proven for
generalized permutohedra in \cite{GPS} to the case of poset cones.

\begin{theorem}
If $\p$ is a poset on the set $[n]$, the antipode $S$ acts on the
quasisymmetric enumerator function $F_q(C(\mathsf{P}))$ by
$$S(F_q(C(\mathsf{P}))=(-1)^n\sum_{\mathcal{G}\in\,\mathfrak{F}(\p)}f(C(\mathsf{P}/\mathcal{G}),-q)M_{\mathcal{G}^{op}}.$$
\end{theorem}

\noindent\textit{Proof.} Since $\mathcal{F}\in\mathfrak{F}(\p)$
and $\mathcal{G}\preceq\mathcal{F}$ implies
$\mathcal{G}\in\mathfrak{F}(\p)$ we have
\begin{equation*}\begin{split}\label{pet}
S(F_q(C(\mathsf{P})))&=\sum_{\mathcal{F}\in\,\mathfrak{F}(\p)}q^{\mathrm{rk}_\p(\mathcal{F})}S(M_{\mathcal{F}})=
\sum_{\mathcal{F}\in\,\mathfrak{F}(\p)}q^{\mathrm{rk}_\p(\mathcal{F})}(-1)^{|\mathcal{F}|+1}\sum_{\mathcal{G}\preceq\mathcal{F}^{op}}M_\mathcal{G}\\
&=\sum_{\mathcal{G}\in\,\mathfrak{F}(\p^{op})}M_{\mathcal{G}}
\sum_{\substack{\mathcal{F}\in\,\mathfrak{F}(\p)
\\\mathcal{G}\,\preceq\mathcal{F}^{op}}}q^{\mathrm{rk}_\p(\mathcal{F})}(-1)^{|\mathcal{F}|+1}.
%&=(-1)^n\sum_{\mathcal{G}\in\,\mathfrak{F}(\p^{op})}M_\mathcal{G}\sum_{\substack{\mathcal{F}\in\,
%\mathfrak{F}(\p)\\\mathcal{G}\,\preceq\mathcal{F}^{op}}}(-q)^{\mathrm{rk}_\p(\mathcal{F})}(-1)^{|\mathcal{F}|+1+n+\mathrm{rk}_\p(\mathcal{F})}
\end{split}\end{equation*}
By equivalences $\mathcal{G}\preceq\mathcal{F}^{op}$ if and only
if $\mathcal{G}^{op}\preceq\mathcal{F}$ and
$\mathcal{G}\in\mathfrak{F}(\p^{op})$ if and only if
$\mathcal{G}^{op}\in\mathfrak{F}(\p)$, the last equality becomes
\begin{equation*}\begin{split}
S(F_q(C(\mathsf{P}))
&=(-1)^n\sum_{\mathcal{G}^{op}\in\,\mathfrak{F}(\p)}M_\mathcal{G}\sum_{\substack{\mathcal{F}\in\,\mathfrak{F}(\p)\\\mathcal{G}^{op}\,
\preceq\mathcal{F}}}(-q)^{\mathrm{rk}_\p(\mathcal{F})}(-1)^{|\mathcal{F}|+1+n+\mathrm{rk}_\p(\mathcal{F})}\\
&=(-1)^n\sum_{\mathcal{G}\in\,\mathfrak{F}(\p)}M_{\mathcal{G}^{op}}\sum_{\substack{\mathcal{F}\in\,\mathfrak{F}(\p)\\\mathcal{G}\,
\preceq\mathcal{F}}}(-q)^{\mathrm{rk}_\p(\mathcal{F})}(-1)^{|\mathcal{F}|+1+n+\mathrm{rk}_\p(\mathcal{F})}.
\end{split}\end{equation*}
By the proof of Theorem \ref{jakobitno}, we have that
$\p/\mathcal{G}$ is positive subposet of $\p$ and $\mathcal{G}$ is
a normal flag to the face $C(\p/\mathcal{G})$. By the same
argument, the faces of $C(\p/\mathcal{G})$ are of the form
$C(\p/\mathcal{F})$ for flags of ideals
$\mathcal{F}\in\mathfrak{F}(\p)$ which satisfy
$\mathcal{G}\preceq\mathcal{F}$. Therefore, according to the
identities (\ref{sssss}) and $\mathrm{rk}_\p(\mathcal{F})=\dim
C(\p/\mathcal{F})$, we obtain
$$f(C(\p/\mathcal{G}),q)=\displaystyle\sum_{\substack{\mathcal{F}\in\,\mathfrak{F}(\p)\\\mathcal{G}\preceq\,\mathcal{F}}}
(-1)^{n+1+|\mathcal{F}|+\mathrm{rk}_\p(\mathcal{F})}q^{\mathrm{rk}_\p(\mathcal{F})}.$$\qed

\section{Function $F(\p)$}

In this section we show that for a well labelled poset the
enumerator function $F_q(C(\mathsf{P}))$ specializes at $q=0$ to
the enumerator of $\mathsf{P}-$partitions. For a poset $\p$ on the
set $[n]$ let
\begin{equation*}
F(\p)=F_0(C(\mathsf{P}))=\sum\limits_{\substack{\mathcal{F}\in\,\mathfrak{F}(\p)\\\p/\mathcal{F}\text{
discrete}}}M_\mathcal{F}.
\end{equation*}
The sum is over flags of ideals of $\p$ such that $\p/\mathcal{F}$
is a \emph{discrete poset} (poset with no covering relations).

\begin{proposition}
The principal specialization  $\bold{ps}^1$ of $F(\mathsf{P})$ at
$m=-1$ results in
\begin{equation*}
\bold{ps}^1(F(\p))(-1)=(-1)^{n-1}.
\end{equation*}
\end{proposition}
\noindent\textit{Proof.} The poset cone $C(\mathsf{P})$ has a
unique vertex, so the required equation follows from Proposition
\ref{f-f-f}.\qed

Let $\p_1\ast\p_2$ be the series composition of posets $\p_1,
\p_2$ which is a poset on the disjoint union of elements of
$\p_1,\p_2$ with the order relation defined by
$i\leq_{\p_1\ast\p_2}j$ if and only if $i\leq_{\p_1}j$ or
$i\leq_{\p_2}j$ or $i\in\p_1, j\in\p_2$.

\begin{proposition}\label{ser_com}
If $\p=\p_1\ast \p_2\ast\cdots\ast \p_n$ then
\begin{equation*}
F(\p)=F(\p_1)\circ F(\p_2)\circ\cdots\circ F(\p_n),
\end{equation*}
where $\circ$ is the concatenation product described in Example
\ref{bipart}.
\end{proposition}
\noindent\textit{Proof.} The proof follows from the fact that
$F(\p)$ is the sum over all flags of ideals $\mathcal{F}$ such
that $\p/\mathcal{F}$ is a discrete poset. Such a flag is the
consecutive composition of flags of the same type corresponding to
components $\p_1,\ldots,\p_n$. \qed

\begin{example}
A poset $\mathsf{K}_{m,n}$ can be expressed as the series
composition of discrete posets $\mathsf{d}(m)\ast \mathsf{d}(n)$
on $[m],[n]$. By Proposition \ref{ser_com}, we obtain
$F(\mathsf{K}_{m,n})=F(\mathsf{d}(m))\circ
F(\mathsf{d}(n))=M^m_{(1)}\circ M^n_{(1)}.$
\end{example}

\begin{proposition}\label{posledicica}
If $\mathrm{Max}(\p)$ is the set of maximal elements of a
connected poset $\mathsf{P}$ on the set $[n]$ then
$$F(\p)=\sum_{\emptyset\neq A\subseteq\mathrm{Max}(\p)}(F(\p|_{[n]\setminus A}))_{|A|}.$$
Particulary, if $\mathrm{Max}(\p)=\{v\}$ then
$F(\p)=\left(F(\p|_{[n]\setminus\{v\}})\right)_1.$
\end{proposition}
\noindent\textit{Proof.} For a flag of ideals
$\mathcal{F}:\emptyset =F_0\subset F_1\subset\cdots\subset
F_m=[n]$ let $F_\mathrm{end}=F_m\setminus F_{m-1}.$ If
$\p/\mathcal{F}$ is a discrete poset, we have that
$\p|_{F_\mathrm{end}}$ is discrete too, so
$F_\mathrm{end}\subseteq\mathrm{Max}(\p).$ Therefore
\begin{equation*}\begin{split}
F(\p)=\sum\limits_{\substack{\mathcal{F}\in\,\mathfrak{F}(\p)\\\p/\mathcal{F}\text{
discrete}}}M_\mathcal{F}=\sum_{\emptyset\neq
A\subseteq\mathrm{Max}(\p)}
\sum\limits_{\substack{\mathcal{F}\in\,\mathfrak{F}(\p)\\\p/\mathcal{F}\text{discrete}\\F_\mathrm{end}=A}}M_\mathcal{F}
=\sum_{\emptyset\neq
A\subseteq\mathrm{Max}(\p)}(F(\p|_{[n]\setminus A}))_{|A|}.
\end{split}\end{equation*}\qed

\begin{remark}
The direct calculation shows that $F(\p)$ distinguishes all
non-isomorphic posets on $[n]$, for $n\in\{1,2,3,4,5,6\}.$
\end{remark}

The weighted enumerator function invariant $F_q$ contains more
information about posets than its specialization at $q=0$. The
following example of posets with the same quasisymmetric invariant
$F(\mathsf{P})$ is borrowed from \cite[Example 4.9]{MW}.

\begin{figure}[h!]
\begin{center}
\includegraphics[width=60mm]{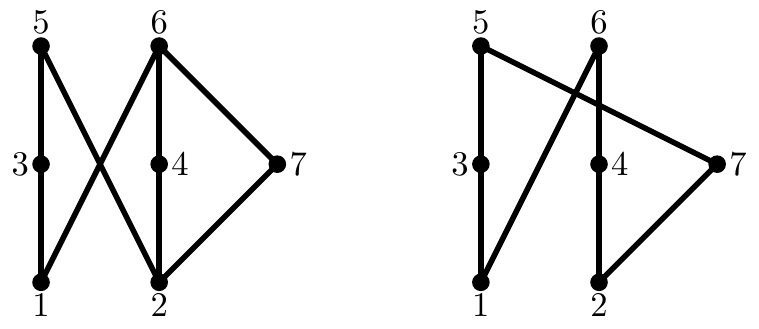}
\caption{Posets with the same quasisymmetric invariant
$F(\mathsf{P})$}
\end{center}
\end{figure}

\begin{example}
Let $\p_1$ and $\p_2$ be posets on Figure 3. The direct
calculation gives
\begin{equation*}\begin{split}
F(\mathsf{\p}_1)=F(\mathsf{\p}_2)&=M_{(2,3,2)}+2M_{(1,1,3,2)}+2M_{(2,3,1,1)}\\
&+M_{(1,3,2,1)}+M_{(1,2,3,1)}+M_{(2,1,3,1)}+M_{(1,3,1,2)}\\
&+3 M_{(1,2,2,2)}+3 M_{(2,1,2,2)}+3 M_{(2,2,1,2)}+3 M_{(2,2,2,1)}\\
&+2 M_{(1,2,2,1,1)}+6 M_{(2,1,1,1,2)}+6 M_{(1,2,2,1,1)}\\
&+7 M_{(2,1,1,2,1)}+7 M_{(1,2,1,1,2)}+8 M_{(1,1,2,1,2)}+8 M_{(2,1,2,1,1)}\\
&+8 M_{(1,2,1,2,1)}+8 M_{(1,1,2,2,1)}+9 M_{(1,1,1,2,2)}+9 M_{(2,2,1,1,1)}\\
&+20 M_{(1,1,1,1,1,2)}+23 M_{(1,1,1,1,2,1)}+24 M_{(1,1,1,2,1,1)}\\
&+24 M_{(1,1,2,1,1,1,)}+23 M_{(1,2,1,1,1,1)}+20 M_{(2,1,1,1,1,1)}\\
&+4 M_{(1,1,3,1,1)}+3 M_{(1,3,1,1,1)}+3 M_{(1,1,1,3,1)}+66
M_{(1,1,1,1,1,1,1)}.
\end{split}\end{equation*}
Consider flags
$\emptyset\subset\{i\}\subset[6]\setminus\{j\}\subset[6],$ for
$i\in\{1,2\}$ and $j\in\{5,6\}$ which are the only flags of ideals
of the type $(1,5,1)$ of posets $\p_1$ and $\p_2$. The
coefficients by $M_{(1,5,1)}$ in $F_q(C(\p_i)),i=1,2$ are
$$\zeta_{(1,5,1)}(\p_1)=2q^4+q^3+q^2\;\;\text{and}\;\;\zeta_{(1,5,1)}(\p_2)=q^4+3q^3,$$
which shows that $F_q(C(\p_1))\neq F_q(C(\p_2))$.
\end{example}

The following theorem gives a geometric interpretation of the
enumerator $F_\p(\bold{x})$ of $\mathsf{P}-$partitions for a well
labelled poset $\mathsf{P}$. Recall from Definition
\ref{definicija} that
\[F_\p(\bold{x})=\sum_{f\in\mathcal{A}(\p)}x_{f(1)}x_{f(2)}\cdots
x_{f(n)},\] where $f\in\mathcal{A}(\p)$ if and only if $i<_\p j$
implies $f(i)<f(j),$ for all $i,j\in\p$.

\begin{proposition}\label{stav_poslednji}
For a well labelled poset $\p$ holds
$$\mathcal{A}(\p)\;=\;\{\omega\in C^\circ_\mathcal{F}\;:\;\mathcal{F}\in\mathfrak{F}(\p)\text{ and }\p/\mathcal{F}\text{ is a discrete poset}\}.$$
\end{proposition}
\noindent\textit{Proof.} We have $\omega\in C^\circ_\mathcal{F}$
for some $\mathcal{F}\in\mathfrak{F}(\p)$ such that
$\p/\mathcal{F}$ is a discrete poset if and only if $\omega^\ast$
is maximized uniquely at the vertex of $C(\p).$ By Proposition
\ref{stav_1} and definition of $\p$-partitions of well labeled
posets, the both sets are characterized by the same condition.
\qed

\begin{theorem}\label{druga}
For a well labelled poset $\mathsf{P}$ holds
$F(\mathsf{P})=F_{\mathsf{P}}(\bold{x}).$
\end{theorem}
\noindent\textit{Proof.} By (\ref{jedinicica}) and Proposition
\ref{stav_poslednji}

\[F(\p)=\sum_{\substack{\mathcal{F}\in\mathfrak{F}(\p)\\\p/\mathcal{F}\text{
discrete}}}M_\mathcal{F}=\sum_{\substack{\omega\,\in\,\mathbb{Z}^n_+\cap\,
C^\circ_\mathcal{F}\\\mathcal{F}\in\mathfrak{F}(\p)\\\p/\mathcal{F}\text{
discrete}}}x_{\omega_1}x_{\omega_2}\cdots
x_{\omega_n}=\sum_{\omega\in\mathcal{A}(\p)}x_{\omega_1}x_{\omega_2}\cdots
x_{\omega_n},\] what is exactly the expression for
$F_\p(\bold{x})$. \qed

\bibliographystyle{amsplain}

\end{document}